\documentclass[12pt]{article}
\usepackage[dvipdfm]{graphicx}          
\usepackage{tikz}

 \usepackage{latexsym}
 \usepackage{amssymb}
 \usepackage{color}
 \usepackage{graphicx}
\usepackage{epstopdf}
\usepackage{threeparttable}
\usepackage{multirow}
\usepackage{amsmath}    
\usepackage{bm}                     
 \newtheorem{Theorem}{Theorem}
\newtheorem{Definition}{Definition}
\newtheorem{Proposition}{Proposition}

\newtheorem{Corollary}{Corollary}

\newcommand{\T}{{\cal T}}

\newcommand{\A}{{\cal A}}
\newcommand{\B}{{\cal B}}

\newcommand{\I}{{\cal I}}

\newcommand{\x}{{\bf x}}

\newcommand{\0}{{\bf 0}}

\newcommand{\qed}{\nobreak \ifvmode \relax \else
      \ifdim\lastskip<1.5em \hskip-\lastskip
      \hskip1.5em plus0em minus0.5em \fi \nobreak
      \vrule height0.75em width0.5em depth0.25em\fi}

\def \ep{\hbox{ }\hfill$\Box$}

\begin{document}
\title{Three Dimensional Strongly Symmetric Circulant Tensors}

\author{Liqun Qi\footnote{Department of Applied Mathematics,
    The Hong Kong Polytechnic University, Hung Hom, Kowloon, Hong Kong.
    E-mail: maqilq@polyu.edu.hk (L. Qi).
    This author's work was partially supported by
    the Hong Kong Research Grant Council (Grant No. PolyU 502111, 501212, 501913 and 15302114).}
    \quad Qun Wang \footnote{Department of Applied Mathematics,
    The Hong Kong Polytechnic University, Hung Hom, Kowloon, Hong Kong.
    Email: wangqun876@gmail.com (Q. Wang).} \quad  Yannan Chen\footnote{School of Mathematics and Statistics,
    Zhengzhou University, Zhengzhou, China. E-mail: ynchen@zzu.edu.cn (Y. Chen)
    This author's work was supported by the
    National Natural Science Foundation of China (Grant No. 11401539) and
    the Development Foundation for Excellent Youth Scholars of
    Zhengzhou University (Grant No. 1421315070).}}

\date{\today} \maketitle

\begin{abstract}
\noindent  
  In this paper, we give a necessary and sufficient condition for
  an even order three dimensional strongly symmetric circulant tensor to be positive semi-definite.
  In some cases, we show that this condition is also sufficient for this tensor to be sum-of-squares.
  Numerical tests indicate that this is also true in the other cases.

\noindent {\bf Key words:}\hspace{2mm} H-eigenvalue,
  strongly symmetric tensor, circulant tensor,
  sum of squares, positive semi-definiteness.

\noindent {\bf AMS subject classifications (2010):}\hspace{2mm} 15A18; 15A69
  \vspace{3mm}

\end{abstract}

\section{Introduction}
\hspace{4mm}

Let $\A = (a_{i_1\cdots i_m})$ be an even order $n$ dimensional real symmetric tensor, i.e., $m=2k$, $i_1,\cdots, i_m = 1, \cdots, n$, $a_{i_1\cdots i_m}$ are real and invariant under any index permutation.   Then $\A$ is corresponding to a homogeneous polynomial $f(\x)$ for $\x \in\Re^n$, defined by
$$f(\x) = \sum_{i_1,\cdots, i_m=1}^n a_{i_1\cdots i_m}x_{i_1}\cdots x_{i_m}.$$
If $f(\x) \ge 0$ for any $\x \in \Re^n$, then $\A$ is called a positive semi-definite (PSD) tensor, $f$ is called a PSD or nonnegative polynomial.   To identify a PSD symmetric tensor or a PSD homogeneous polynomial is an important problem in theory and practice \cite{AP, Hi, Qi}.   It is NP-hard to identify a general even order symmetric tensor is PSD or not.   Recently, it was discovered that several easily checkable classes of special even order symmetric tensors are PSD.    These include even order symmetric diagonally dominated tensors \cite{Qi}, even order symmetric
B$_0$ tensors \cite{QS}, even order Hilbert tensors \cite{SQ}, even order symmetric M tensors \cite{ZQZ}, even order symmetric double B$_0$ tensors \cite{LL},
even order symmetric strong H tensors \cite{LWZZL, KSB}, even order strong Hankel tensors \cite{Qi15}, even order positive Cauchy tensors \cite{CQ}, etc.

Hankel tensors arise from signal processing and data fitting \cite{BB, DQW, PDV, Qi15, Xu}.  They are symmetric tensors.   It is known that even order strong Hankel tensors are PSD.    A Hankel tensor is called a strong Hankel tensor if its generating vector generates a PSD Hankel matrix.    This condition is checkable.  However, there is still no checkable condition for verifying a general even order Hankel tensor is PSD or not.
In \cite{LQX}, it was proved that an even order strong Hankel tensor is a sum-of-squares (SOS) tensor.   It was also discovered there that there are PSD Hankel tensors, which are not strong Hankel tensors, but still SOS tensors.   Thus, an open question is raised in \cite{LQX}: are all PSD Hankel tensors SOS tensors? If the answer to this question is ``yes'', then the problem for determining a given even order Hankel tensor is PSD or not can be solved by solving a semi-definite linear programming  problem \cite{LQX, Las, Lau}.

It is not easy to answer the question raised in \cite{LQX}.   For general homogeneous polynomials, this problem was first studied by Hilbert \cite{Hi}.
In 1888, Hilbert \cite{Hi} proved that only in the following three cases, a PSD homogeneous polynomial, of degree $m$ and $n$ variables, definitely is an SOS polynomial: 1) $m = 2$; 2) $n = 2$; 3) $m=4$ and $n=3$, where $m$ is the degree of the polynomial and $n$ is the number of variables.    Hilbert proved that in all the other possible combinations of $n$ and even $m$, there are PSD non-SOS (PNS) homogeneous polynomials.    However, Hilbert did not give an explicit example of PNS homogeneous polynomials.   The first explicit example of PNS homogeneous polynomials was given by Motzkin \cite{Mo} in 1967.   More examples of PNS homogeneous polynomials can be found in \cite{CL, Re}.

According to Hilbert \cite{Hi, Re}, two cases with low values of $m$ and $n$, in which there are PNS homogeneous polynomials, are the case that $m=6$ and $n=3$, and the case that $m=n=4$.     In \cite{LQW} and \cite{CQW}, sixth order three dimensional Hankel tensors and fourth order four dimensional Hankel tensors were studied respectively.    No PNS Hankel tensors were found there.  However, even for such low values of $m$ and $n$, a strict proof to show that PNS Hankel tensors do not exist seems still very difficult.    In these two cases, the Hankel tensors have thirteen independent entries of their generating vectors.   Even assuming that their generating vectors are symmetric \cite{CQW}, there are still seven independent entries.    The situation is still complicated enough to make a strict proof \cite{CQW}.

Are there other special classes of even order symmetric structured tensors, there are no easily checkable conditions for their positive semi-definiteness, but it is possible that they are PNS-free?    A good candidate for such a problem is the class of even order strongly symmetric circulant tensors.

Strongly symmetric tensors were introduced in \cite{QXX}.   An $m$th order $n$ dimensional tensor $\A = (a_{i_1\cdots i_m})$ is called a strongly symmetric tensor if
$$a_{i_1\cdots i_m} \equiv a_{j_1\cdots j_m}$$
as long as $\{ i_1, \cdots, i_m \} = \{ j_1, \cdots, j_m \}$.   Note that a symmetric matrix is a strongly symmetric tensor.   Hence, strongly symmetric
tensors are also extensions of symmetric matrices.   Some good properties of symmetric matrices may be inherited by strongly symmetric tensors, not symmetric tensors.

Circulant tensors have applications in stochastic process and spectral hypergraph theory \cite{CzQ}.   An $m$th order $n$ dimensional tensor $\A = (a_{i_1\cdots i_m})$ is called a circulant tensor if
$$a_{i_1\cdots i_m} \equiv a_{j_1\cdots j_m}$$
as long as $j_l = i_l+1, mod(n)$ for $l=1, \cdots, m$.

Hence, in this paper, we consider even order three dimensional strongly symmetric circulant tensors.

We find that a general three dimensional strongly symmetric circulant tensor has only three independent
entries: the diagonal entry $d$, the off-diagonal entry $u$, which has two different indices,  and the off-diagonal entry $c$, which has
three different indices.   Let $\A = (a_{i_1\cdots i_m})$ be an $m$th order three dimensional strongly symmetric circulant tensor.  Here $m$ can be even or odd.    We denote
$$a_S \equiv a_{i_1\cdots i_m}$$
if $S = \{i_1, \cdots, i_m \}$.
Since $\A$ is strongly symmetric, it has at most seven independent entries for $S = \{1\}, \{2\}, \{ 3\}, \{1, 2\}, \{1, 3\}, \{2, 3\}$ and $\{1, 2, 3\}$.
Since $\A$ is also circulant, we have
$$a_{\{1\}} = a_{\{2\}} = a_{\{3\}}, \ {\rm and}\ a_{\{1, 2\}} = a_{\{2, 3\}} = a_{\{1, 3\}}.$$
Let $d = a_{\{1\}} = a_{\{2\}} = a_{\{3\}}$, $u= a_{\{1, 2\}} = a_{\{2, 3\}} = a_{\{1, 3\}}$ and $c = a_{\{1, 2, 3\}}$.   Then we see that $d$ is the diagonal entry of $\A$: $d = a_{1\cdots 1} = a_{2\cdots 2} = a_{3\cdots 3}$, and an $m$th order three dimensional strongly symmetric circulant tensor has only three independent entries $d, u$ and $c$.  Thus, we may denote a general three dimensional strongly symmetric circulant tensor $\A = \A(m, d, u, c)$, where $m$ is its order.   When the context is clear, we only use $\A$ to denote it.

In our discussion, we need the concept of H-eigenvalues of symmetric
tensors, which was introduced in \cite{Qi} and is closely related to
positive semi-definiteness of even order symmetric tensors. In the
next section, we introduce H-eigenvalues and
discuss their relations with positive semi-definiteness of even
order symmetric tensors.

Now, let $m=2k$ be even.   In Section 3,
we show that there are two one-variable functions $M_c(u)$ and $N_c(u)$,
such that $M_{c}(u) \ge N_{c}(u) \ge 0$, $\A$ is SOS if and only if $d \ge M_c(u)$,
and $\A$ is PSD if and only if $d \ge N_c(u)$.
If $M_c(u) = N_c(u)$, then
three dimensional strongly symmetric PNS circulant tensors do not exist
for such $u$ and $c$.   We show that if $u, c \le 0$ or $u = c > 0$, then $M_c(u) = N_c(u)$.   Explicit formulas for $M_c(u) = N_c(u)$ are given there in these cases.  Thus, it is PNS-free for such $u$ and $c$.

Note that $\A$ is PSD or SOS if and only if $\alpha \A$ is PSD or SOS, respectively.   Thus, we only need to consider three cases that $c =0$, $c=1$ and $c=-1$.

In Section 4, we discuss the case that $c=0$.  In this case, for $u>0$, we have $M_0(u) = uM_0(1)$ and $N_0(u) = uN_0(1)$.   We show that $-N_0(1)$ is the smallest H-eigenvalue of $\A(m, 0, 1, 0)$.   Numerical tests show that $M_0(1) = N_0(1)$ for $m = 6, 8, 10, 12$ and $14$.

In Section 5, we study the case that $c = -1$.  We show that there is a $u_0 > 0$ such that if $u \le u_0$, $N_{-1}(u)$ is linear and the explicit formula of $N_{-1}(u)$ can be given, and if $u > u_0$, $N_{-1}(u)$ is the smallest H-eigenvalue of a tensor with $u$ as a parameter.    Numerical tests show that for $u > 0$, we still have $M_{-1}(u) = N_{-1}(u)$ for $m = 6, 8, 10$ and $12$.

In Section 6, we study the case that $c = 1$.  We show that there is a $v_0 < 0$ such that if $u \le v_0$, $N_1(u)$ is linear and the explicit formula of $N_1(u)$ can be given, and if $u > v_0$, $N_1(u)$ is the smallest H-eigenvalue of a tensor with $u$ as a parameter.    Numerical tests show that for $u \not = 1$, we still have $M_1(u) = N_1(u)$ for $m = 6, 8, 10$ and $12$.

Some final remarks are made in Section 7.

\section{H-eigenvalues}
\hspace{4mm}

H-eigenvalues of symmetric tensors were introduced in \cite{Qi}.   They are closely related to positive semi-definiteness of even order symmetric tensors.   Let $\T = (t_{i_1\cdots i_m})$ be an $m$th order $n$ dimensional real symmetric tensor and $\x \in \Re^n$.   Then $\T \x^{m-1}$ is a vector in $\Re^n$, with its $i$th component defined as
$$\left(\T \x^{m-1}\right)_i = \sum_{i_2, \cdots, i_m=1}^n t_{ii_2\cdots i_m}x_{i_2}\cdots x_{i_m}.$$
If there is $\x \in \Re^n, \x \not = \0$ and $\lambda \in \Re$ such that for $i=1, \cdots, n$,
$$\left(\T \x^{m-1}\right)_i = \lambda x_i^{m-1},$$
then $\lambda$ is called an H-eigenvalue of $\T$ and $\x$ is called its associated H-eigenvector.   When $m$ is even, H-eigenvalues always exist.   $\T$ is PSD if and only if its smallest H-eigenvalue is nonnegative \cite{Qi}.   From now on, we denote the smallest H-eigenvalue of $\A(m, d, u, c)$ as
$\lambda_{\min}(m, d, u, c)$.

\section{Functions $M_c(u)$ and $N_c(u)$}
\hspace{4mm}

 In this section and the next three sections, we assume
that $n=3$ and $m=2k$ is even.   Let $\A$ be an $m$th order three
dimensional strongly symmetric circulant tensor.    Then we may
write $f_c(\x) \equiv f(\x) = \A \x^m$ as
\begin{eqnarray}
f_c(\x) & =  & d(x_1^m + x_2^m + x_3^m) + u\sum_{p=1}^{m-1} \left({m \atop p}\right)(x_1^{m-p}x_2^p + x_1^{m-p}x_3^p + x_2^{m-p}x_3^p) \nonumber \\ && +c\sum_{p=1}^{m-2}\sum_{q=1}^{m-p-1}\left({m \atop p}\right)\left({m-p \atop q}\right)x_1^{m-p-q}x_2^px_3^q. \label{e1}
\end{eqnarray}

We now establish two functions $M_c(u)$ and $N_c(u)$, in the following theorem.  Recall that for an $m$th order $n$ dimensional tensor $\A = a_{i_1\cdots \i_m}$, its $i$th off-diagonal entry absolute value sum is defined as
$$r_i = \sum_{i_2, \cdots, i_m=1}^n |a_{ii_2\cdots i_m}| - | a_{ii\cdots i}|.$$
If $a_{i\cdots i} \ge r_i$ for $i=1, \cdots, n$, then $\A$ is called diagonally dominated, and all of its H-eigenvalues \cite{Qi} are nonnegative.  If furthermore $\A$ is even order and symmetric, then $\A$ is PSD \cite{Qi} and SOS \cite{CLQ}.

\begin{Theorem}  \label{t1}
Let $\A$ be an $m$th order three dimensional strongly symmetric circulant tensor.
Then, there are two convex functions $M_c(u) \geq N_c(u) \ge 0$
such that $\A$ is SOS if and only if $d \ge M_c(u)$,
and $\A$ is PSD if and only if $d \ge N_c(u)$.  If
$M_c(u) = N_c(u)$, then $m$th order three dimensional PNS strongly symmetric circulant tensors do not exist for such $u$ and $c$.
Furthermore, we have
\begin{equation} \label{e2}
    M_c(u) \le |u|(2^m-2) +|c|(3^{m-1}-2^m+1).
\end{equation}
\end{Theorem}


\noindent
{\bf Proof}   Since $\A$ is a circulant tensor, then it has the same off-diagonal entry absolute value sum for different rows, i.e., $r_1 = r_2 = r_3$.   By (\ref{e1}), this row sum is equal to the right hand side of (\ref{e2}).  Thus, if $d$ is greater than or equal to this value, $\A$ is diagonally dominated and thus PSD and SOS.   This shows the existence of $N_c(u)$, $M_c(u)$ and (\ref{e2}). As the set of PSD tensors and the set of SOS tensors are convex \cite{LQX}, $M_c(u)$ and $N_c(u)$ are convex.  Since a necessary condition for an even order circulant tensor to be PSD is that its diagonal entry to be nonnegative \cite{CzQ}, we have $N_c(u) \ge 0$ for all $u$ and $c$.

By definition, we have $N_c(u) \le M_c(u)$.  Clearly, if
$M_c(u) = N_c(u)$, then $m$th order three dimensional PNS strongly symmetric circulant tensors do not exist for such $u$ and $c$.   The theorem is proved.
\ep

As discussed in the introduction, for the PNS-free problem, we only need to consider three values of $c$: $c = 0, 1$ and $-1$.

If $M_c(u) = N_c(u)$, $u$ is called a {\bf PNS-free point} for $c$.

For the convenience, we present formally three ingredients used in theoretical proofs of PNS-free points.
If a point $u$ enjoys these ingredients, it is PNS-free.

\begin{Definition}
Suppose that $n=3$ and $m=2k$ is even.   Suppose that there is a number $M$ such that $\A$ is SOS if $d = M$, and
a nonzero vector $\bar \x \in \Re^3$ such that
$f_c^*(\bar \x) = 0$, where $f_c^*(\x) \equiv f_c(\x)$ with $d = M$.
Then we call $M$ the {\bf critical value} of $\A$ at $u$, the SOS decomposition $f_c^*(\x)$ the {\bf critical SOS decomposition} of $\A$ at $P$, and $\bar \x$ the {\bf critical minimizer} of $\A$ at $u$.
\end{Definition}

\begin{Theorem} \label{t2}
Let $u \in \Re$.   Then $u$ is PNS-free for $c$ if $\A$ has a critical value $M$, a critical SOS decomposition $f_c^*(\x)$ and a critical minimizer $\bar{\x}$ at $u$.
\end{Theorem}

\noindent
{\bf Proof}  Suppose that $\A$ has a critical value $M$, a critical SOS decomposition $f_c^*(\x)$ and a critical minimizer $\bar{\x}$ at $u$.   Then we have $M \ge M_c(u)$ by the definition of $M_c(u)$.   If $d < M$, then
$$f_c(\bar \x) = (d-M)(\bar x_1^m+\bar x_2^m +\bar x_3^m) + f_c^*(\bar \x) < 0.$$
This implies that $N_c(u) \ge M$ by the definition of $N_c(u)$.  But $N_c(u) \le M_c(u)$.  Thus, $M_c(u) = N_c(u) = M$, i.e., $u$ is PNS-free for $c$.
\ep

\begin{Corollary} \label{c1}
If $u, c \le 0$, then
\begin{equation} \label{e3}
M_c(u) = N_c(u) = -u(2^m-2) -c(3^{m-1}-2^m+1).
\end{equation}
Thus, it is PNS-free for such $u$ and $c$.
\end{Corollary}

\noindent
{\bf Proof}  Suppose that $u, c \le 0$.   Let $M$ be the value of the right hand side of (\ref{e2}), and $\bar \x = (1, 1, 1)^\top$.   If $d = M$, then
$f_c(\x) = f_c^*(\x)$ has an SOS decomposition as $\A$ is an even order diagonally dominated symmetric tensor \cite{CLQ}.   We also see that $f_c^*(\bar \x) = 0$.  The result follows.
\ep

\begin{Corollary} \label{c2}
If $u = c > 0$, then
$$M_c(u) = N_c(u) = u = c.$$
Thus, it is PNS-free for such $u$ and $c$.
\end{Corollary}

\noindent
{\bf Proof}  Suppose that $u = c > 0$.   Let $M = u = c$, and $\bar \x = (2, -1, -1)^\top$.   If $d = M$, then
$f_c(\x) = f_c^*(\x) = (x_1+x_2+x_2)^m$ has an SOS decomposition.   We also see that $f_c^*(\bar \x) = 0$.  The result follows.
\ep

\begin{Corollary} \label{c3}
If $u > 0$, then
$$M_0(u) = uM_0(1)$$
and
$$N_0(u) = uN_0(1).$$
Hence, for $c=0$, it is PNS-free if and only if $M_0(1) = N_0(1)$.
\end{Corollary}

\noindent
{\bf Proof}  Suppose that $u > 0$ and $d \ge uM_0(1)$.   By (\ref{e1}), we have
\begin{eqnarray*}
f_0(\x) & = & (d-uM_0(1))(x_1^m + x_2^m + x_3^m)\\
 &&+ u\left(M_0(1)(x_1^m + x_2^m + x_3^m) + \sum_{p=1}^{m-1} \left({m \atop p}\right)(x_1^{m-p}x_2^p + x_1^{m-p}x_3^p + x_2^{m-p}x_3^p)\right).
\end{eqnarray*}
We see that $f_0(\x)$ is SOS.  Hence, $M_0(u) = uM_0(1)$.    Similarly, we may prove that $N_0(u) = uN_0(1)$.   By these and Corollary \ref{c1}, we have the last conclusion.
\ep


\section{$c=0$}
\hspace{4mm}

If $u \le 0$, by Corollary \ref{c1}, we have $M_0(u) = N_0(u)= -u(2^m -2)$.
If $u > 0$, by Corollary \ref{c3}, we have $M_0(u) = uM_0(1)$ and $N_0(u) = uN_0(1)$.
We only need to consider the case that $u=1$.

\begin{Proposition} \label{p1}
We have that $N_0(1) = -\lambda_{\min}(m, 0, 1, 0)$.
\end{Proposition}
\noindent
{\bf Proof}  By \cite{Qi}, $\A(m, d, 1, 0)$ is PSD if and only if $\lambda_{\min}(m, d, 1, 0) \ge 0$.   By the structure of circulant tensors,
$\lambda_{\min}(m, d, 1, 0) = d + \lambda_{\min}(m, 0, 1, 0)$.   Thus, $\A(m, d, 1, 0)$ is PSD if and only if $d \ge - \lambda_{\min}(m, 0, 1, 0)$.
By the definition of $N_c(u)$, we have $N_0(1) = -\lambda_{\min}(m, 0, 1, 0)$.
\ep

For $m = 6, 8, 10, 12$ and $14$,  we compute $M_0(1)$ and $N_0(1)$ by using Matlab (YALMIP, GloptiPloy and SeDuMi) software and Maple \cite{HLL, Lo, L2, St}, respectively.
We find for such $m$, $M_0(1) = N_0(1)$.  The results are displayed in Table 1.

\begin{table}[!htb]
  \centering
  \begin{tabular}{c|cc}
    \hline\hline
    $m$ & $M_0(1)$ & $N_0(1)$ \\
    \hline
6  &1.737348471173345 & 1.737348471777547 \\
8  &1.882980354978972 & 1.882980356780414\\
10 &1.947977161918168 & 1.947977172341075\\
12 &1.976878006619490 & 1.976878047128592\\
14 &1.989722829997529 & 1.989723542124766\\
    \hline\hline
  \end{tabular}
  \caption{The values of $M_0(1)$ and $N_0(1)$.}\label{Table-1}
\end{table}

\section{$c=-1$}
\hspace{4mm}

If $u \le 0$, then Corollary \ref{c1} indicates that $M_{-1}(u) = N_{-1}(u) = -u(2^m-2) +(3^{m-1}-2^m+1)$.  We now discuss the case that $u > 0$.

In this section and the next section, we denote that
$\B =\A(m, 3^{m-1}-2^m+1, 0, -1)$ and $\T = \A(m, 2^m-2,-1,0)$.
Then, $\B$ and $\T$ are obviously diagonally dominated. Hence, they are PSD and SOS \cite{CLQ}.
And all of their H-eigenvalues are nonnegative.

\begin{Theorem}  \label{t3}
Let $$ \varphi(u)\equiv \lambda_{\min}(\B-u\T), $$
where $\lambda_{\min}(\cdot)$ denotes the smallest H-eigenvalue.
Then, $\varphi(u) \le 0$.   If $\varphi(u) = 0$, then we have
\begin{equation} \label{e4}
N_{-1}(u) = 3^{m-1}-2^m+1 -u(2^m-2).
\end{equation}
If $\varphi(u) < 0$, then we have
\begin{equation} \label{e5}
N_{-1}(u) = 3^{m-1}-2^m+1 -u(2^m-2) -\lambda_{\min}(m, 3^{m-1}-2^m+1 -u(2^m-2), u, -1).
\end{equation}
Furthermore, the set $C = \{ u : \varphi(u) = 0 \}$ is a nonempty closed convex ray $(-\infty, u_0]$ for some $u_0 \ge 0$.
\end{Theorem}
\noindent
{\bf Proof}
Let $\bar x = (1, 1, 1)^\top$.   Then $\B \bar \x^m = 0$ and $\T \bar \x^m = 0$.    Thus, $(\B-u\T)\bar \x^m = 0$ for any $u$.   By \cite{Qi}, we see that
$$ \varphi(u)\equiv \lambda_{\min}(\B-u\T) \le 0.$$

If $\varphi(u) = 0$, let $d = 3^{m-1}-2^m+1 -u(2^m-2)$.  Then $\A(m, d, u, -1) = \B-u\T$.  We have $\lambda_{\min}(m, d, u, -1) = 0$.  This implies (\ref{e4}).

If $\varphi(u)<0$, then, because
\begin{equation*}
    f_{-1}(\x)=(d-(3^{m-1}-2^m+1)+u(2^m-2))\I\x^m+\B\x^m-u\T\x^m \geq 0,
\end{equation*}
we have
\begin{eqnarray*}
  N_{-1}(u) &=& \inf\{d~:~\lambda_{\min}((d-(3^{m-1}-2^m+1)+u(2^m-2))\I+\B-u\T)\geq 0\} \\
   &=& (3^{m-1}-2^m+1)-u(2^m-2)-\lambda_{\min}(\B-u\T) \\
   &=& 3^{m-1}-2^m+1 -u(2^m-2) -\lambda_{\min}(m, 3^{m-1}-2^m+1 -u(2^m-2), u, -1).
\end{eqnarray*}
 We have (\ref{e5}).

By Corollary \ref{c1}, $C$ is nonempty and $u \in C$ as long as $u \le 0$.   By Theorem \ref{t1}, $N_{-1}(u)$ is a convex function.  By this, (\ref{e4}) and (\ref{e5}), $C$ is convex.  Since $\lambda_{\min}$ is a continuous function \cite{Qi}, $C$ is closed.   Since $u \in C$ as long as $u \le 0$, $C$ is a ray, with the form $(-\infty, u_0]$ for some $u_0 \ge 0$.
\ep

\begin{Corollary} \label{c4}
  Let $u_0\equiv \max \{\hat{u}~:~\varphi(\hat{u})=0\}.$   Then $u_0$ is well-defined and
  $u_0\geq 0$.    Furthermore, for $u \le u_0$, we have (\ref{e4}),
  and for $u > u_0$, we have (\ref{e5}).
\end{Corollary}

\begin{Proposition} \label{p2}
If $M_{-1}(u_0)= N_{-1}(u_0) = 3^{m-1}-2^m+1 -u_0(2^m-2)$, then for $u \le u_0$, we have
$M_{-1}(u) = N_{-1}(u) = 3^{m-1}-2^m+1 -u(2^m-2)$.
\end{Proposition}

\noindent
{\bf Proof}  By Theorem \ref{t1}, $M_{-1}(u)$ is convex.   By Corollary 1, $M_{-1}(u) = 3^{m-1}-2^m+1-u(2^m-2)$ for $u \le 0$.    Since $u_0 \ge 0$, the
conclusion follows.
\ep

\begin{Proposition}
  Suppose  $u_0= \max \{\hat{u}~:~\varphi(\hat{u})=0\}.$ Then, we
  have
  \begin{equation} \label{e6}
  0 \leq u_0 \leq \bar{u}_0(m) \equiv \frac{3^{m-1}+1}{2^m} -1.
  \end{equation}
\end{Proposition}
\noindent {\bf Proof}
  Since $\B$ is PSD and has a H-eigenvalue $0$, we have
  $\varphi(0)=0$ and $u_0\geq 0$.

  On the other hand, we consider the case $u>\bar{u}_0$.
  Let $\x_0 = (1,1,-3)^{\top}$. We have
  $$ (\B-\bar{u}_0\T)\x_0^m = 0 \qquad\text{ and }\qquad \T\x_0^m = 2^m(3^m-1). $$
  Then,
  \begin{equation*}
    (\B-u\T)\x_0^m = (\B-\bar{u}_0\T)\x_0^m - (u-\bar{u}_0)\T\x_0^m
      = -(u-\bar{u}_0)2^m(3^m-1) <0.
  \end{equation*}
  Hence, we have $\varphi(u)=\lambda_{\min}(\B-u\T)<0$ when
  $u>\bar{u}_0$. Therefore, $u_0 \leq \bar{u}_0$.
\ep

For $m = 6, 8, 10, 12$ and $14$, we find that $\B-\bar u_0\T$ is PSD.   This shows that for such $m$, $\varphi(\bar u_0) = 0$, i.e.,
\begin{equation} \label{e6.1}
  u_0 = \bar{u}_0(m) \equiv \frac{3^{m-1}+1}{2^m} -1.
  \end{equation}
It remains a further research topic to show that $\B-\bar u_0\T$ is PSD for all even $m$ with $m \ge 16$.   If this is true, then
(\ref{e6.1}) is true for all even $m$ with $m \ge 6$.


\begin{table}[!htb]
  \centering
  \begin{tabular}{c|cc}
    \hline\hline
    $u$ & $M_{-1}(u)$ & $N_{-1}(u)$ \\
    \hline
0.1 &173.799999999899& 173.8\\
2 &55.9999999995172& 56\\
${45 \over 16}$ &5.62499991033116& 5.625\\
5 &9.42544641511067& 9.4254465011842588\\
10 &18.1121860822789& 18.112186280892696\\
40 &70.2326321651344& 70.232638183914150\\
300 &521.943237017699&521.94324013633004 \\
    \hline\hline
  \end{tabular}
  \caption{The values of $M_{-1}(u)$ and $N_{-1}(u)$ for $m=6$.}\label{Table-2}
\end{table}

\begin{table}[!htb]
  \centering
  \begin{tabular}{c|cc}
    \hline\hline
    $u$ & $M_{-1}(u)$ & $N_{-1}(u)$ \\
    \hline
1 &1677.99999992219& 1678\\
3 &1169.99999999356& 1170\\
${483 \over 64}$ &15.0937478786308& 15.09375\\
10 &19.7129359300341& 19.7129361640501\\
40 &76.2023466001335& 76.2023468071730\\
140 &264.500365037152& 264.500382469583\\
300 &565.777184078832& 565.777239551091 \\
    \hline\hline
  \end{tabular}
  \caption{The values of $M_{-1}(u)$ and $N_{-1}(u)$ for $m=8$.}\label{Table-3}
\end{table}

\begin{table}[!htb]
  \centering
  \begin{tabular}{c|cc}
    \hline\hline
    $u$ & $M_{-1}(u)$ & $N_{-1}(u)$ \\
    \hline
1 &17637.9999999549& 17638\\
10 &8439.99999783081& 8440\\
${4665 \over 256}$ &36.4452603485520& 36.4453125\\
20 &39.9075358817909& 39.9075375522954\\
40 &78.8670625326286& 78.8670809985488\\
140 &273.664775923815& 273.664798232238\\
300 &585.341085323688& 585.341145806726 \\
    \hline\hline
  \end{tabular}
  \caption{The values of $M_{-1}(u)$ and $N_{-1}(u)$ for $m=10$.}\label{Table-4}
\end{table}

\begin{table}[!htb]
  \centering
  \begin{tabular}{c|cc}
    \hline\hline
    $u$ & $M_{-1}(u)$ & $N_{-1}(u)$ \\
    \hline
1 &168957.999979042& 168958\\
20 &91171.9999996683& 91172\\
${43263 \over 1024}$ &84.4971787852022& 84.498046875 \\
60 &119.589505579120& 119.589562756497 \\
100 &198.664532858285& 198.664684641639 \\
140 &277.739708996851& 277.739806526784 \\
300 &594.040191670531& 594.040294067366 \\
    \hline\hline
  \end{tabular}
  \caption{The values of $M_{-1}(u)$ and $N_{-1}(u)$ for $m=12$.}\label{Table-5}
\end{table}

In Tables 2-5, the values of $M_{-1}(u)$ and $N_{-1}(u)$ for $m=6,8,10,12$ and  $u=0.1,2,{45 \over 16},5,10,40,300$, $u=1,3,{483 \over 64},10,40,140,300$, $u=1,10,{4665 \over 256},20,40,140,300$ and $u=1,20,{43263 \over 1024},60,100,140,300$ are reported, respectively.   We find for such $m$ and $u$, $M_{-1}(u) = N_{-1}(u)$.

\newpage

\section{$c=1$}
\hspace{4mm}

Corollary \ref{c2} indicates that $M_1(1) = N_1(1) = 1$.   Hence, we only need to consider the case that $u \not = 1$.   Let $\B$ and $\T$ be the same as in the last section.    We have the following
theorem.

\begin{Theorem} \label{t4}
Let $$ \psi(u)\equiv \lambda_{\min}(-u\T-\B). $$
Then, $\psi(u) \le 0$.   If $\psi(u) = 0$, then we have
\begin{equation} \label{e7}
N_1(u) = -(3^{m-1}-2^m+1) -u(2^m-2).
\end{equation}
If $\psi(u) < 0$, then we have
\begin{equation} \label{e8}
N_1(u) = -(3^{m-1}-2^m+1) -u(2^m-2) -\lambda_{\min}(m, -(3^{m-1}-2^m+1) -u(2^m-2), u, 1).
\end{equation}
Furthermore, if the set $C = \{ u : \psi(u) = 0 \}$ is nonempty, then it is a closed convex ray $(-\infty, v_0]$ for some $v_0 < 0$.
\end{Theorem}

\noindent
{\bf Proof}   The proof of this theorem is similar to the proof of Theorem \ref{t3}.   However, we cannot apply Corollary \ref{c1} here.
If $C$ is non empty,
we may show that $C$ is closed and convex as in the proof of Theorem \ref{t3}.

If there is a $\hat{u}\leq0$ such that $\lambda_{\min}(-\hat{u}\T-\B)=0$, for $u\leq\hat{u}\leq0$, we have
\begin{eqnarray*}
  \psi(u) &=& \lambda_{\min}(-u\T-\B) \\
    &=& \lambda_{\min}(-\hat{u}\T-\B+(-u+\hat{u})\T) \\
    &\geq& \lambda_{\min}(-\hat{u}\T-\B) \\
    &=& 0.
\end{eqnarray*}
Hence, $\psi(u)=0$ for all $u\leq\hat{u}\leq0$.
So if $C$ is not empty, it is a ray with the form  $(-\infty, v_0]$ for some $v_0$.   Clearly, $\psi(0) \equiv \lambda_{\min}(-\B) < 0$ as $\B$ is PSD and not a zero tensor.   Hence, $v_0 < 0$.

The other parts of the proof are similar to the proof of Theorem \ref{t3}.
\ep

\begin{Corollary} \label{c5}
   If there is one point $\hat u$ such that $\psi(\hat u) = 0$, let $v_0\equiv \max \{\hat{u}~:~\psi(\hat{u})=0\}$.
 Then for $u \le v_0$, we have (\ref{e7}),
  and for $u \ge v_0$, we have (\ref{e8}).
\end{Corollary}

We also have the following proposition.

\begin{Proposition} \label{p3}
If $M_{1}(v_0)= N_{1}(v_0) = -(3^{m-1}-2^m+1) -v_0(2^m-2)$, then for $u \le v_0$, we have
$M_{1}(u) = N_{1}(u) = -(3^{m-1}-2^m+1) -u(2^m-2)$.
\end{Proposition}

\noindent
{\bf Proof}  Suppose that $M_{1}(v_0)= N_{1}(v_0) = -(3^{m-1}-2^m+1) -v_0(2^m-2)$.  By (\ref{e3}), if $u \le v_0$ and $d=-(3^{m-1}-2^m+1) -u(2^m-2)$,
we have
$$f_{1}^*(\x) = -\bar ug_1(\x) + g_2(\x),$$
where
$$g_1(\x) = \A(m, 2^m-2, -1, 0),$$
$$g_2(\x) = \A(m, -(3^{m-1}-2^m+1)-v_0(2^m-2), v_0, 1)$$
and
$$\bar u \equiv u-v_0 \le 0.$$
We see that $g_2(\x)$ is equal to the critical SOS decomposition of $\A$ at $c= 1$ and $u = v_0$, and $g_1(\x)$ is equal to the critical SOS decomposition of $\A$ at $c=0$ and $u=-1$.   Hence both $g_1(\x)$ and $g_2(\x)$ are SOS polynomials.
This implies that $f_{1}^*(\x)$ is an SOS polynomial.   Let $\bar \x = (1, 1, 1)^\top$, we see that $f_{1}^*(\bar \x) = 0$.   Then the conclusion follows from Theorem \ref{t2}.
\ep

\begin{Proposition}
  Suppose that $C$ is not empty.  Let $v_0 = \max \{\hat{u}~:~\psi(\hat{u})=0\}.$ Then, we
  have
  \begin{equation} \label{e9}
  v_0 \leq \bar{v}_0(m) \equiv 1-\frac{3^{m-1}}{2^{m-1}+1}.
  \end{equation}
\end{Proposition}
\noindent {\bf Proof}
  Let $\x_0 = (1,1,-\frac{1}{2})^{\top}$. We have
  $$ (-\B-\bar{v}_0\T)\x_0^m = 0 \qquad\text{ and }\qquad \T\x_0^m = 2^m+1-2^{1-m}. $$
  Then,
  \begin{equation*}
    (-\B-u\T)\x_0^m = (-\B-\bar{v}_0\T)\x_0^m - (u-\bar{v}_0)\T\x_0^m
      = -(u-\bar{v}_0)(2^m+1-2^{1-m}) <0.
  \end{equation*}
  Hence, we have $\psi(u)=\lambda_{\min}(-\B-u\T)<0$ when
  $u>\bar{v}_0$. Therefore, $v_0 \leq \bar{v}_0$.
\ep

Similar to the discussion of $u_0$, for $m = 6, 8, 10, 12$ and $14$, we find that $-\bar v_0\T-\B$ is PSD.   This shows that for such $m$, $\psi(\bar v_0) = 0$, i.e.,
\begin{equation} \label{e9.1}
 v_0 = \bar{v}_0(m) \equiv 1-\frac{3^{m-1}}{2^{m-1}+1}.
  \end{equation}
This also shows that $C$ is not empty for such $m$.
It remains a further research topic to show that $-\bar v_0\T-\B$ is PSD for all even $m$ with $m \ge 16$.   If this is true, then
(\ref{e9.1}) is true for all even $m$ with $m \ge 6$.

In Tables 6-9, the values of $M_1(u)$ and $N_1(u)$ for $m = 6, 8, 10, 12$ and $u=-40,-10,-{70\over11},\\-5,-1,10,40$, $u=-40,-20,-{686\over43},-10,-1,10,40$, $u=-100,-40,-{710\over19},-20,-1,10,40$ and $u=-140,-100,-{58366\over683},-60,-1,10,40$ are reported, respectively.   We find for such $m$ and $u$, $M_1(u) = N_1(u)$.

\begin{table}[!htb]
  \centering
  \begin{tabular}{c|cc}
    \hline\hline
    $u$ & $M_1(u)$ & $N_1(u)$ \\
    \hline
-40 &2299.99999993444& 2300\\
-10 &439.999999987168& 440\\
$-{70\over11}$ &214.545454213196& 214.54545454545454\\
-5 &173.991050854352& 173.99105151869704\\
-1 &55.8846973214056& 55.884697712412670\\
10 &16.6347897201042&16.634789948247836 \\
40 &68.7552353830704& 68.755241242186800\\
    \hline\hline
  \end{tabular}
  \caption{The values of $M_1(u)$ and $N_1(u)$ for $m=6$.}\label{Table-6}
\end{table}

\begin{table}[!htb]
  \centering
  \begin{tabular}{c|cc}
    \hline\hline
    $u$ & $M_1(u)$ & $N_1(u)$ \\
    \hline
-40 &8228.00000000754& 8228 \\
-20 &3147.99999997053& 3148 \\
$-{686 \over 43}$ &2120.18604151092& 2120.18604651163 \\
-10 &1371.80748977461& 1371.80749709544 \\
-1 &243.740078469126& 243.740080311110 \\
10 &17.9466697015668& 17.9466711544215 \\
40 &74.4360734714431& 74.4360817805826\\
    \hline\hline
  \end{tabular}
  \caption{The values of $M_1(u)$ and $N_1(u)$ for $m=8$.}\label{Table-7}
\end{table}

\begin{table}[!htb]
  \centering
  \begin{tabular}{c|cc}
    \hline\hline
    $u$ & $M_1(u)$ & $N_1(u)$ \\
    \hline
-100 &83539.9999994888& 83540\\
-40 &22219.9999995661& 22220\\
$-{710 \over 19}$ &19530.5255392283& 19530.5263157895 \\
-20 &10678.1156381743& 10678.1156702343 \\
-1 &1004.40451207948& 1004.40454284172 \\
10 &18.5317674915799& 18.5317776218259 \\
40 &76.9710868541002& 76.9710927899669\\
    \hline\hline
  \end{tabular}
  \caption{The values of $M_1(u)$ and $N_1(u)$ for $m=10$.}\label{Table-8}
\end{table}

\begin{table}[!htb]
  \centering
  \begin{tabular}{c|cc}
    \hline\hline
    $u$ & $M_1(u)$ & $N_1(u)$ \\
    \hline
-140 &400107.999992756& 400108\\
-100 &236347.999998551& 236348\\
$-{58366\over683}$ &176802.173593347& 176802.170881802\\
-60 &124727.840916646& 124727.840144917\\
-1 &4063.38103939314& 4063.38106552746 \\
10 &18.7918976770375& 18.7919005425937 \\
40 &78.0982265029468& 78.0982419563963 \\
    \hline\hline
  \end{tabular}
  \caption{The values of $M_1(u)$ and $N_1(u)$ for $m=12$.}\label{Table-9}
\end{table}

\newpage

\section{Final Remarks}
\hspace{4mm}

In Proposition \ref{p1}, Theorems \ref{t3} and \ref{t4}, we give a necessary and sufficient condition for
an even order three dimensional strongly symmetric circulant tensor to be positive semi-definite.
For $u, c \le 0$ and $u=c > 0$, we show that this condition is also sufficient for this tensor to be sum-of-squares.
Numerical tests indicate that this is also true in the other cases.

How can $\B-\bar u_0\T$ and $-\bar v_0\T-\B$ be shown to be PSD for all even $m \ge 6$?  If these are true, then
(\ref{e6.1}) and (\ref{e9.1}) are true for all even $m \ge 6$.

Finally, more efforts are needed to prove that this problem is PNS-free eventually.

\vspace{3mm}
\hspace{4mm}
\hspace{4mm}

%
%


\end{document}